\documentclass[12]{amsart} 
\usepackage{amssymb,latexsym}
\usepackage{array}
\usepackage{enumerate}

\usepackage{amsthm}
\begin{document}
\title{Some results on the topology of real Bott towers} 
\author{ Raisa Dsouza and V.Uma}
\address{Department of Mathematics, IIT Madras,
  Chennai, India} 
\email{raisadsouza1989@gmail.com; vuma@iitm.ac.in} 
\date{}

\begin{abstract}  
  The main aim of this article is to study the topology of real Bott
  towers as special and interesting examples of real toric
  varieties. We first give a presentation of the fundamental group of
  a real Bott tower and show that the fundamental group is abelian if
  and only if the real Bott tower is a product of circles. We further
  prove that the fundamental group of a real Bott tower is always
  solvable and it is nilpotent if and only if it is abelian. We then
  describe the cohomology ring of a real Bott tower and also give
  recursive formulae for the Steifel Whitney classes. We derive
  combinatorial characterization for orientability of these manifolds
  and further give a combinatorial formula for the $(n-1)$th Steifel
  Whitney class. In particular, we show that if a Bott tower is
  orientable then the $(n-1)$th Steifel Whitney class must also
  vanish. Moreover, by deriving a combinatorial formula for the second
  Steifel-Whitney class we give a necessary and sufficient condition
  for the Bott tower to admit a spin structure. We finally prove the
  vanishing of all the Steifel-Whitney numbers and hence establish
  that these manifolds are null-cobordant.
\end{abstract}
\maketitle 

\thispagestyle{empty} 

\footnote{{\bf keywords:} Real Bott tower, cohomology ring,
fundamental group, Steifel-Whitney classes}

\def\theequation {\arabic{section}.\arabic{equation}}

\newcommand{\be}{\begin{equation}}
\newcommand{\ee}{\end{equation}}
\newcommand{\codim}{\mbox{{\rm codim}$\,$}}
\newcommand{\stab}{\mbox{{\rm stab}$\,$}}
\newcommand{\lr}{\mbox{$\longrightarrow$}}

\newcommand{\ch}{{\mathcal H}}
\newcommand{\cf}{{\mathcal F}}
\newcommand{\cd}{{\mathcal D}}

\newcommand{\ra}{\rightarrow}
\newcommand{\blr}{\Big \longrightarrow}
\newcommand{\da}{\Big \downarrow}
\newcommand{\ua}{\Big \uparrow}
\newcommand{\hra}{\mbox{{$\hookrightarrow$}}}
\newcommand{\rt}{\mbox{\Large{$\rightarrowtail$}}}
\newcommand{\dua}{\begin{array}[t]{c}
\Big\uparrow \\ [-4mm]
\scriptscriptstyle \wedge \end{array}}

\newtheorem{guess}{Theorem}[section]
\newcommand{\bth}{\begin{guess}$\!\!\!${\bf }~}
\newcommand{\eeth}{\end{guess}}
\renewcommand{\bar}{\overline}
\newtheorem{propo}[guess]{Proposition}
\newcommand{\bpropo}{\begin{propo}$\!\!\!${\bf }~}
\newcommand{\epropo}{\end{propo}}

\newtheorem{lema}[guess]{Lemma}
\newcommand{\blem}{\begin{lema}$\!\!\!${\bf }~}
\newcommand{\elem}{\end{lema}}

\newtheorem{defe}[guess]{Definition}
\newcommand{\bdefe}{\begin{defe}$\!\!\!${\bf }~}
\newcommand{\edefe}{\end{defe}}

\newtheorem{coro}[guess]{Corollary}
\newcommand{\bcor}{\begin{coro}$\!\!\!${\bf }~}
\newcommand{\ecor}{\end{coro}}

\newtheorem{rema}[guess]{Remark}
\newcommand{\brem}{\begin{rema}$\!\!\!${\bf }~\rm}
\newcommand{\erem}{\end{rema}}

\newtheorem{exam}[guess]{Example}
\newcommand{\beg}{\begin{exam}$\!\!\!${\bf }~\rm}
\newcommand{\eeg}{\end{exam}}

\newtheorem{notn}[guess]{Notation}
\newcommand{\bnot}{\begin{notn}$\!\!\!${\bf}~\rm}
\newcommand{\enot}{\end{notn}}

\newcommand{\ctext}[1]{\makebox(0,0){#1}}
\setlength{\unitlength}{0.1mm}
\newcommand{\cl}{{\mathcal L}}
\newcommand{\cp}{{\mathcal P}}
\newcommand{\ci}{{\mathcal I}}
\newcommand{\cR}{{\mathcal R}}

\newcommand{\bz}{\mathbb{Z}}
\newcommand{\cs}{{\mathcal s}}
\newcommand{\cv}{{\mathcal V}}
\newcommand{\ce}{{\mathcal E}}
\newcommand{\ck}{{\mathcal K}}
\newcommand{\cz}{{\mathcal Z}}
\newcommand{\cg}{{\mathcal G}}
\newcommand{\bq}{\mathbb{Q}}
\newcommand{\bt}{\mathbb{T}}
\newcommand{\bh}{\mathbb{H}}
\newcommand{\br}{\mathbb{R}}
\newcommand{\wt}{\widetilde}
\newcommand{\im}{{\rm Im}\,}
\newcommand{\bc}{\mathbb{C}}
\newcommand{\bp}{\mathbb{P}}
\newcommand{\spin}{{\rm Spin}\,}
\newcommand{\ds}{\displaystyle}
\newcommand{\tor}{{\rm Tor}\,}
\newcommand{\bff}{{\bf F}}
\newcommand{\bs}{\mathbb{S}}
\def\ns{\mathop{\lr}}
\def\nssup{\mathop{\lr\,sup}}
\def\nsinf{\mathop{\lr\,inf}}
\renewcommand{\phi}{\varphi}
\newcommand{\co}{{\cal O}}
\noindent

\section{Introduction}

Bott towers are iterated fibre bundles with fibre at each stage being
$\bp_\bc^1$. In particular they are smooth projective complex toric
varieties. They were constructed in \cite{gk} by M.Grossberg and
Y.Karshon who show that a Bott-Samelson variety can be deformed to a
Bott tower.  Bott-Samelson manifolds were first constructed in
\cite{bs} to study cohomology of generalized flag
varieties. M.Demazure and D.Hansen used it to obtain
desingularizations of Schubert varieties in generalized flag
varieties. Moreover, the underlying differentiable structure is
preserved under the deformation. Because of their relation with Bott
Samelson manifolds which in turn are related to the Schubert
varieties, along with their amenable structure as iterated
$\bp_\bc^1$-bundles, the Bott towers have been important and
interesting objects of study.

Topological invariants for these manifolds have been studied using
their iterated sphere bundle structure.  The equivariant cohomology of
Bott Samelson manifolds have been studied with applications to the
cohomology of Schubert varieties in \cite{du}.  Also see \cite{w} for
equivariant $K$-theory of Bott towers with application to the
equivariant $K$-theory of flag manifolds.

Indeed, from the viewpoint of toric topology, Bott-towers can be seen
to also have the structure of a quasi-toric manifold \cite{dj} with
the quotient polytope being the $n$-dimensional cube $I^n$ where $n$
is the complex dimension of the Bott tower. The second named author
and P.Sankaran described the topological $K$-ring of the quasi toric
manifolds in \cite{su}, where the topological $K$-ring of Bott towers
and Bott Samelson manifolds have been described as a special
example. Also see \cite{cr}, \cite{br} and \cite{bb} for results on
complex $K$-theory, $KO$-theory as well as the complex cobordism ring
of Bott towers and quasi-toric manifolds.

Recently there has been extensive work on the topology and geometry of
Bott towers viewed as a quasi-toric manifold (see for example
\cite{mp}, \cite{cms}). These works are especially related
to the problem of cohomological rigidity of Bott manifolds or more
generally of quasi-toric manifolds.

There has also been a parallel study on the topology of real Bott
towers. These manifolds are constructed as iterated
$\bp^1_{\br}=\bs^1$-bundles and can be viewed as a special example of
a small cover defined by Davis and Januszkiewicz (see for example
\cite{m1}, \cite{km}, \cite{cmo} ).

The cohomology of a small cover have been described in
\cite{dj}. Indeed in \cite[Theorem 5.12]{dj} Davis and Januszkiewicz
give a presentation for the cohomology ring with $\bz_2$-coefficients
as a quotient of the Stanley-Reisner ring of the simple convex
polytope by certain canonical linear relations. A similar description
for the cohomology ring with $\bz_2$-coefficients of the real part of
a smooth projective complex toric variety is due to Jurkiewicz (see
\cite[Theorem 4.3.1]{j}).

In \cite[Section 6]{dj}, there is also a description of the tangent bundle of a
small cover and a formula for its Steifel-Whitney class in terms of
the generators of the cohomology ring. These generators are in turn
the first Steifel-Whitney classes of certain canonical line bundles on
the small cover.

In \cite{u}, the second named author of the current article described
the fundamental group of a real toric variety associated to any smooth
fan $\Delta$ in $\bz^n$. She in particular gave a presentation of the
fundamental group, a combinatorial criterion on $\Delta$ for the
fundamental group to be abelian along with a criterion for the real
toric variety to be aspherical. The fundamental group of a small cover
has also earlier been described in \cite[Corollary 4.5]{dj} as the
kernel of a natural map from the right-angled Coxeter group associated
to the simple convex polytope $P$ to $\bz_2^n$. The natural map is the
composition of the characteristic map $\lambda$ of the small cover
with the abelianization map of the Coxeter group.

The authors in the current article were interested in the study of
real Bott towers as a nice class of small covers or real toric
varieties, which are characterized by the upper triangular real {\it
  Bott matrix} $C=(c_{i,j})$ with entries in $\bz_2$ (see definition
below).

Our main motivation here is to give a precise and elegant description
of all the above mentioned topological invariants namely, the
fundamental group, cohomology ring and the Stiefel Whitney classes, in
terms of the {\it Bott numbers} $\{c_{i,j}\}$. We exploit the
inductive definition of this class of manifolds and derive precise
topological information about them wherever possible. We further
classify Bott towers satisfying a specific topological property like
orientability or admission of spin structure, by means of certain
algebraic identities on the $c_{i,j}$'s.

We now develop some notations before outlining our main results in the
next section.
\subsection{Notations and Conventions}

In this section we recall the definition of a Bott tower and fix some
notations (see \cite{gk}).

A {\it Bott tower} is a smooth complete complex toric variety which
is constructed iteratively as follows:

Let $Y_1=\bc\bp^1$. Let $L_2$ be a holomorphic complex line bundle on
$\bc\bp^1$.  We then let $Y_2=\bp(\mathbf{1}\bigoplus L_2)$ where
$\mathbf{1}$ is the trivial line bundle on $\bc\bp^1$. Then $Y_2$ is a
$\bc\bp^1$ bundle over $\bc\bp^1$ which is a Hirzebruch surface. We
can iterate this process for $2\leq j\leq n$, where at each step,
$L_j$ is a complex line bundle over $Y_{j-1}$, and the variety
$Y_j=\bp(\mathbf{1}\bigoplus L_j)$ is a $\bc\bp^1$ bundle over
$Y_{j-1}$. The variety $Y_n$ thus obtained after $n$-steps is called
an $n$-step Bott tower.

\bdefe\label{Bottfan}
In fact an $n$-step Bott tower is a smooth complete toric variety of
dimension $n$ whose fan $\Delta$ can be described as follows:

We take a collection of integers $\{c_{i,j}\}, ~1\leq i<j\leq n$. Let
$e_1,e_2,\ldots,e_n$ be the standard basis vectors of $\br^n$. Let
$v_j=e_j$ for $1\leq j\leq n$,
\[v_{n+j}=-e_j+\sum_{k=j+1}^{n}c_{j,k}\cdot e_k\] for $1\leq j\leq n-1$ and
$v_{2n}=-e_n$.  We define the fan $\Delta$ in $\br^n$ consisting of
cones generated by the set of vectors in any subcollection of
$\{v_1,v_2,\ldots,v_n,v_{n+1},\ldots,v_{2n}\}$ which does not contain
both $v_i$ and $v_{n+i}$ for $1\leq i\leq n$.  \edefe

\bdefe\label{Bottqt} We can also view a Bott tower as a quasi-toric
manifold (see \cite{dj}) over the $n$-cube $I^n$ which is a simple
convex polytope of dimension $n$. If we index the $2n$ facets of $I^n$
by $F_1,F_2,\ldots,F_n,F_{n+1},\ldots,F_{2n}$, then the characteristic
function is defined on the collection of facets to $\bz^n$ as follows:
$\lambda(F_j)=e_j$ for $1\leq j\leq n$,
\[\lambda(F_{n+j})=-e_j+\sum_{k=j+1}^nc_{j,k}\cdot e_{j+k}\] for
$1\leq j\leq n-1$ and $\lambda(F_{2n})=-e_n$.  \edefe

\subsubsection{Real Bott tower}
We shall call the real part of the $n$-step complex Bott tower as the
real $n$-step Bott tower. 

In particular, $(Y_2)_{\br}$ is an $\br\bp^1$ bundle over
$(Y_1)_{\br}=\br\bp^1$. Iteratively we construct $(Y_j)_{\br}$ as an
$\br\bp^1$ bundle over $(Y_{j-1})_{\br}$ for $2\leq j\leq n$.  The
real $n$-step Bott tower $(Y_n)_{\br}$ is indeed the real toric
variety associated to the fan $\Delta$ described above (see
\cite[Section 2.4]{j}
and \cite{u}).

\bdefe\label{rbsc} As in the complex case we can also view
$(Y_n)_{\br}$ as a {\it small cover} over the simple convex polytope
$I^n$, where the characteristic map $\lambda$ is defined on the
collection of facets $\cf$ to ${\bz}_2^n$ as follows:
$\lambda(F_j)=e_j$ for $1\leq j\leq n$, 
\[\lambda(F_{n+j})=e_j+\sum_{k=j+1}^nc_{j,k}\cdot e_k\] for
$1\leq j\leq n-1$ and $\lambda(F_{2n})=e_n$. Here $c_{i,j}\in\bz_2$
for $1\leq i<j\leq n$. Thus $(Y_n)_{\br}$ is homeomorphic to the
identification space ${\bz_2^n\times I^n}/\sim$ where
$(t,p)\sim(t^{\prime},p^{\prime})$ if and only if $p=p^{\prime}$ and
$t\cdot (t^{\prime})^{-1}\in G_{F(p)}$. Here
$F(p)=F_1\cap\cdots\cap F_l$ is the unique face of $I$ which contains
$p$ in its relative interior and $G_{F(p)}$ is the rank-$l$ subgroup
of $\bz_2^n$ determined by the span of
$\lambda(F_1),\ldots,\lambda(F_l)$. Now, let $\pi:(Y_n)_{\br}\ra I^n$
denote the second projection which maps $[t,p]\mapsto p$, and let
$M_i:=\pi^{-1}(F_i)$ denote the {\it characteristic submanifold} for
$1\leq i\leq 2n$. (See Section 1 of \cite{dj}).  \edefe

The topological structure of an $n$-step real Bott tower is completely
determined by the simple convex polytope $I^n$ and the data encoded by
the matrix \be\label{BM} C=(c_{i,j})\in M_n(\bz_2)\ee where
$c_{i,i}=1$ and $c_{i,j}=0$ for $i>j$. Note that the $i$th row of $C$
is $\lambda(F_{n+i})\in \bz_2^n$ for $1\leq i\leq n$. We call $C$ the
{\it Bott matrix}. Thus $Y_n=Y(C)$ the real Bott tower associated to
$C$.

The $2$-step real Bott tower is the torus or the Klein bottle
depending on whether $c_{1,2}=0$ or $c_{1,2}=1$.  The $3$-step real
Bott tower is an $\br\bp^1$ bundle over the torus or the Klein bottle
whose topological structure depends on $c_{1,2},c_{1,3}$ and $c_{2,3}$.

\bnot  In this article, since we are mainly interested
in the study of the real Bott tower, for notational simplicity we shall
henceforth denote $(Y_n)_{\br}$ by $Y_n$.\enot

\subsection{Overview of the main results}
In Section 2 we give a description of the fundamental group of the
real Bott tower. In particular, we give a presentation of the
fundamental group in terms of generators and relations in Theorem
\ref{fun} and in Corollary \ref{ab} we prove that $Y_n$ has abelian
fundamental group if and only if $Y_n$ is a product of
circles. Further, in Proposition \ref{sol} we show that the commutator
subgroup of $\pi_1(Y_n)$ is always abelian, so that $\pi_1{(Y_n)}$ is
solvable. In Proposition\ref{nil} we prove that the fundamental group
is nilpotent if and only if it is abelian. By abelianizing the
fundamental group we further determine the first homology with integer
coefficients $H_{1}(Y_n;\bz)$ explicitly. We also conclude by
induction that $Y_n$ is always aspherical. The main tool used in this
section is the presentation of the fundamental group of the real part
of any smooth toric variety in \cite{u}, where combinatorial
characterizations are also given for the fundamental group to be
abelian and the manifold to be aspherical. However, as mentioned
earlier, in the case of the real Bott tower the fundamental group gets
a neater and simpler presentation in terms of the entries of the
characterizing Bott matrix. Moreover, the presentation also enables us
to make further conclusions about the group theoretic properties of
the fundamental group of these special class of real toric
varieties. Further, proofs of all results in this section except that
of Theorem \ref{fun} which uses \cite[Lemma 3.2]{u}, are made
self-contained and specific to the case of the real Bott tower.

In Section 3 we study the cohomology ring of these manifolds with
$\bz_2$-coefficients. The tool here is to apply the description of the
cohomology ring of small covers in \cite{dj} or else that of a real
toric variety in \cite[Section 2.4]{j} to these class of manifolds as a special
case. We describe the cohomology ring in terms of generators and
relations in Theorem \ref{mainiso}, where again the presentation given
in terms of the entries of the Bott matrix becomes simple in this case
due to the iterative structure of these manifolds. 

In Section 4, using the presentation of the cohomology ring and
applying the more general results of \cite{dj} we give explicit
description of the Steifel-Whitney classes of $Y_n$ which are smooth
connected compact manifolds. The description is again in terms of the
entries of the Bott matrix. More explicitly, in Theorem \ref{rf} we
give an inductive formula for total Steifel-Whitney class as well as
the $k$th Steifel-Whitney class of $Y_n$ in terms of those of
$Y_{n-1}$.  

These results are applied to give a combinatorial characterisation for
orientability of $Y_n$ in Theorem \ref{crorient1}, criterion for $Y_n$
to admit a spin structure in Theorem \ref{crspin} and a nice formula
for $w_{n-1}(Y_n)$ in Corollary \ref{n-1th sfclass}. We further show
in Corollary \ref{conseq} that if $Y_n$ is orientable then
$w_{n-1}(Y_n)$ also vanishes. We finally prove the vanishing of all
Steifel-Whitney numbers of $Y_n$, which enables us to conclude in
Theorem \ref{vanishswnos} that the real Bott towers are null
cobordant.

We aim to continue these methods to prove results, for instance, about
the immersion and embedding dimensions and parallelizability of these
manifolds in future work.

{\bf Acknowledgement:} The authors wish to thank Prof. P. Sankaran for
valuable discussions and comments. The first author wishes to thank
University Grants Commission (UGC), India for financial assistance.

\section{The fundamental group}
In this section we shall get a presentation for the fundamental group
of the real Bott tower $Y(\Delta)=Y_{n}$ where $\Delta$ is as in  Definition \ref{Bottfan}.

For obtaining a presentation for $\pi_1(Y_n)$ we apply \cite[Proposition
3.1]{u}. We shall first fix some notations which are slightly modified
from \cite{u} suitable to our setting.  Recall that we have an exact
sequence \be\label{esf} 0\ra\pi_1(Y_n)\ra W(\Delta)\ra \bz_2^n\ra 0\ee
where $W(\Delta)$ is the right angled Coxeter group associated to
$\Delta$ with the following presentation:
\[W(\Delta)=\big{\langle} s_j\mid s_j^2,~ 1\leq j\leq 2n ~~\mbox{and}
~ (s_is_j)^2~ \mbox{for all}~ 1\leq i< j\leq 2n ~\mbox{with}~ j \neq
i+n \big{\rangle} .\]
The last arrow in the above exact sequence is obtained by composing
the natural abelianization map from $W(\Delta)$ to $\bz^{2n}$ with the
characteristic map $\lambda$ from $\bz^{2n}$ to $\bz^n$.

Let
$\alpha_j:=s_js_{j-n}s_{j-n+1}^{c_{j-n, j-n+1}}\cdots
s^{c_{j-n,j-n+k}}_{j-n+k} \cdots s_{n}^{c_{j-n, n}}$
for all $n+1\leq j\leq 2n-1$ and $1\leq k\leq 2n-j$ and
$\alpha_{2n}:=s_{2n}s_n$.

Let $b^j_{i}:=c_{j-n,i}$ for $j-n+1\leq i\leq n$, $b^j_{j-n}=1$ and
$b^j_k=0$ for $k<j-n$, for every $n+1\leq j\leq 2n$. Thus \be\label{bvector} b^j
=(b^j_i)_{i=1,\ldots,n}\ee denotes the $(j-n)$th row vector of the
Bott matrix $C$.

For $\epsilon=(\epsilon_1,\epsilon_2,\ldots \epsilon_n)\in \bz_2^n$
where $\epsilon_i\in\{0,1\}$, let
\[t_{\epsilon} := s_1^{\epsilon_1}\cdots s_n^{\epsilon_n}\]
in $W(\Delta)$.

Now, from the relations in $W(\Delta)$ it follows
that: \be\label{1} t_{\epsilon}\alpha_j
  t_{\epsilon}=\left\{\begin{array}{ll} \alpha_j & \mbox{if~~}
  \epsilon_{j-n}=0\\\\\alpha_j^{-1}& \mbox{if~~}
  \epsilon_{j-n}=1 \end{array}\right. \ee

Since $b^j_{j-n}=1$ and $b^j_{k-n}=0$ for $k-n<j-n$ it follows that:
\be\label{2} t_{b^j}\alpha_jt_{b^j}=\alpha_{j}^{-1}\ee for all
$n+1\leq j\leq 2n$ and \be\label{3} t_{b^j}\alpha_kt_{b^j}=\alpha_{k} \ee for
$n+1\leq k<j \leq 2n$. Moreover, since $b^{j}_{k-n}=c_{j-n,k-n}$ we
have:
\be\label{4} t_{b^j}\alpha_k t_{b^j}=\left\{\begin{array}{ll} 
\alpha_k & \mbox{if~~}c_{j-n,k-n} =0\\\\\alpha_k^{-1}& \mbox{if~~}
c_{j-n,k-n}=1 \end{array}\right.
\ee
for all $n+1\leq j<k \leq 2n$.

\bth\label{fun} We have a presentation for $\pi_1(Y_n)=\langle S ~|~R \rangle$ where
\be\label{gens} S=\{\alpha_j ~:~n+1\leq j\leq 2n\}\ee and  \be\label{csrels} R=\{x_{p,q}
~;~n+1\leq p<q \leq 2n\}\ee where
\begin{displaymath}
  x_{p,q} =
  \left\{\begin{array}{ll} \alpha_p\alpha_q\alpha_p^{-1}
           \alpha_q^{-1} & \mbox{if}~~ c_{p-n,q-n}=0\\\\
           \alpha_p\alpha^{-1}_q\alpha_p^{-1}\alpha_q^{-1}& \mbox{\it if}~~
                                                            c_{p-n,q-n}=1
\end{array}\right\} .\end{displaymath}\eeth
{\bf Proof:} Since $v_1=e_1,\ldots,v_n=e_n$ form a basis of
$N\otimes_{\bz}\bz_2$ and also pairwise form cones in $\Delta$, we can
apply  \cite[Lemma 3.2]{u}. We can relate the above notations with
those in the proof of  \cite[Lemma 3.2]{u} as follows:
\[y_{j,\epsilon}=t_{\epsilon}\cdot \alpha_j\cdot t_{\epsilon}\]
for $n+1\leq j\leq 2n$ and $\epsilon\in \bz_2^n$. Thus
$y_{j,\epsilon}=\alpha_j^{\pm1}$ by (\ref{1}). Let $B^p:=\epsilon+b^p$
and $B^q:=\epsilon+b^q$ in $\bz_2^n$, where $b^p$ and $b^q$ are as
defined in (\ref{bvector}). Then it can be seen that
\[y_{p,\epsilon}\cdot y_{q,B^{q}}\cdot y_{p,B^{p}+B^{q}}\cdot
y_{q,B^{q}}=t_\epsilon \cdot\alpha_p\cdot t_{b^p}\alpha_q t_{b^p}\cdot
t_{b^q+b^p}\alpha_pt_{b^q+b^p}\cdot t_{b^q}\alpha_qt_{b^q}\cdot
t_{\epsilon}\]
whenever $n+1\leq p,q \leq 2n$ Moreover, by (\ref{1}) (\ref{2}),
(\ref{3}) and (\ref{4}) we can further see that if $p<q$ and
$c_{p-n,q-n}=0$ then \[\begin{array}{lll}\alpha_p\cdot
  t_{b^p}\alpha_q t_{b^p}\cdot t_{b^q+b^p}\alpha_pt_{b^q+b^p}\cdot
  t_{b^q}\alpha_qt_{b^q}&=(\alpha_q\cdot t_{b^q}\alpha_p t_{b^q}\cdot
                          t_{b^q+b^p}\alpha_qt_{b^q+b^p}\cdot t_{b^p}\alpha_pt_{b^p} )^{-1}\\
                        &=\alpha_p\alpha_q\alpha_p^{-1}\alpha_q^{-1}\end{array}\]
                      and if $p<q$ and $c_{p-n,q-n}=1$ then
\[\begin{array}{lll}\alpha_p\cdot t_{b^p}\alpha_q t_{b^p}\cdot
    t_{b^q+b^p}\alpha_pt_{b^q+b^p}\cdot
    t_{b^q}\alpha_qt_{b^q}&=(\alpha_q\cdot t_{b^q}\alpha_p t_{b^q}\cdot
                            t_{b^q+b^p}\alpha_qt_{b^q+b^p}\cdot t_{b^p}\alpha_pt_{b^p} )^{-1}\\
                          &=\alpha_p\alpha^{-1}_q\alpha_p^{-1}\alpha_q^{-1}.\end{array}\] Also
                        it can be seen that
                        when $p=q$ we simply get a trivial relation.
                        Thus by \cite[Lemma 3.2]{u} it follows that that $\pi_1(Y_n)$ has a
                        presentation with generators \be S=\{\alpha_j ~:~n+1\leq j\leq
                        2n~\}\ee and the relations
\[R'=\{x^{\epsilon}_{p,q} ~;~n+1\leq p<q \leq 2n\}\] where
\be\label{fullrels} x^{\epsilon}_{p,q} = \left\{\begin{array}{ll}
    t_{\epsilon}\alpha_p\alpha_q\alpha_p^{-1}
                                                  \alpha_q^{-1}t_{\epsilon} & \mbox{if}~~ c_{p-n,q-n}=0\\\\
                                                  t_{\epsilon}\alpha_p\alpha^{-1}_q\alpha_p^{-1}\alpha_q^{-1}
                                                  t_{\epsilon}&
                                                                \mbox{\it
                                                                if}~~
                                                                c_{p-n,q-n}=1
         \end{array}\right\}\ee for $\epsilon\in\bz_2^n$. Moreover, it
       can also be seen that each $x^{\epsilon}_{p,q}$ is conjugate to
       either $x_{p,q}$ or $x^{-1}_{p,q}$ by an element of the free
       group on $S$. For instance, in the case when $\epsilon_{p-n}=0$
       and $\epsilon_{q-n}=1$ and $c_{p-n,q-n}=1$ we have:
       \[\begin{array}{lll} x^{\epsilon}_{p,q} &= t_{\epsilon}\cdot\alpha_p\alpha^{-1}_q\alpha_p^{-1}\alpha_q^{-1}
                                                 \cdot
                                                 t_{\epsilon}\\ &=t_{\epsilon}\alpha_pt_{\epsilon}\cdot
                                                                  t_{\epsilon}\alpha^{-1}_qt_{\epsilon}\cdot
                                                                  t_{\epsilon}\alpha_p^{-1}
                                                                  t_{\epsilon}
                                                                  \cdot
                                                                  t_{\epsilon}
                                                                  \alpha_q^{-1}
                                                                  t_{\epsilon}\\
                                               &=\alpha_p\alpha_q\alpha_p^{-1}\alpha_q\\   &=
                                                                                             (\alpha^{-1}_q)\cdot
                                                                                             (\alpha_p\alpha^{-1}_q\alpha_p^{-1}\alpha_q^{-1})^{-1}
                                                                                             \cdot(\alpha_q)\\
                                               &=\alpha^{-1}_q\cdot
                                                 x^{-1}_{p,q}\cdot
                                                 \alpha_q.\end{array}\]
                                               The other cases can be
                                               proved similarly (see Appendix).
                                               Thus it follows that $R$ given by (\ref{csrels}) is a complete
                                               set of relations for $\pi_{1}(Y_n)$. Hence the theorem. $\Box$

                                               \bcor\label{ab} The
                                               group $\pi_1(Y_n)$ is
                                               abelian if and only if
                                               $Y_n$ is a product of
                                               $n$-copies of $\bp^1$.
                                               \ecor {\bf Proof:} From
                                               the above theorem we
                                               see that $\pi_1(Y_n)$
                                               is abelian if and only
                                               if $c_{p-n,q-n}=0$ for
                                               all
                                               $n+1\leq p < q\leq 2n$.
                                               Hence the
                                               corollary. $\Box$

\blem\label{aspher} $Y_n$ is an aspherical manifold.  \elem {\bf
  Proof: } Note that $Y_n$ is an $\bs^1$-bundle over $Y_{n-1}$ so by
the long exact homotopy sequence we get:
\be\label{lehs}\ldots\ra\pi_{2}(Y_{n-1})\ra
\pi_1(\bs^1)\ra\pi_{1}(Y_n)\ra\pi_{1}(Y_{n-1})\ra 0\ldots.\ee Since
$\pi_{k}(Y_{1})=0$ for $k\geq 2$, the proof follows by induction on
$n$. $\Box$

\bcor The group $\pi_1(Y_n)$ is always torsion free.  \ecor {\bf
  Proof:} This follows readily since $Y_n$ is an aspherical
manifold. $\Box$

\brem Lemma \ref{aspher} alternately follows from \cite[Theorem
1.2]{u}, since it can easily be seen that the fan $\Delta$ associated
to $Y_n$ is {\it flag like}. More precisely, we cannot find three edge
vectors $\{v_i,v_j,v_k\}$ which pairwise form a cone in $\Delta$ but
together do not form a cone in $\Delta$.  Also, Corollary \ref{ab} can
alternately be derived from the combinatorial criterion on $\Delta$
for the fundamental group to be abelian in \cite[Theorem 5.1]{u}.\erem

\subsection{Further group theoretic properties of $\pi_1(Y_n)$}

\bpropo\label{sol} The commutator subgroup $[\pi_{1}(Y_n),\pi_{1}(Y_n)]$ is
abelian. In particular, $\pi_1(Y_n)$ is a solvable group.  \epropo
{\bf Proof:} For every $v_i$, $1\leq i\leq n$ (resp.
$n+1\leq i\leq 2n$) there exists a unique $v_{i+n}$ (resp. $v_{i-n}$)
such that $v_i,v_{i+n}$ (resp. $v_i,v_{i-n}$) do not form a cone in
$\Delta$. Thus by \cite[Lemma 4.1]{u}, $[W,W]$ is abelian. Further,
since $\pi_{1}(Y_n)$ is a subgroup of $W$ (see \ref{esf}) it follows
that $[\pi_{1}(Y_n),\pi_1(Y_n)]$ is a subgroup of $[W,W]$ and is hence
abelian. Now, $1\leq[\pi_1,\pi_1]\leq\pi_1(Y_n)$ gives an abelain
tower for $\pi_1(Y_n)$ so that $\pi_{1}(Y_n)$ is solvable. $\Box$

Let $\bar{\alpha_j}$ denote the image of $\alpha_j$ under the
canonical abelianization homomorphism
\[\pi_{1}(Y_n)\ra H_1(Y_n;\bz)\simeq
\pi_{1}(Y_n)/[\pi_1(Y_n),\pi_1(Y_n)].\]
We then have the following description of $H_{1}(Y_n;\bz)$:
\bpropo\label{presfirsthomology} The group $H_1(Y_n;\bz)$ has a
presentation with generators \be\label{gens}\langle \bar{\alpha_j} :
n+1\leq j\leq 2n \rangle\ee and relations
\be\label{rel1}\bar{\alpha_p}\cdot\bar{\alpha_q}\cdot\bar{\alpha_p}^{-1}\cdot\bar{\alpha_q}^{-1}\ee
for $n+1\leq p,q\leq 2n$ and \be\label{rel2} \bar{\alpha_q}^2\ee for
those $n+1\leq q\leq 2n$ for which there exists a $p<q$ such that
$c_{p-n,q-n}=1$.  Thus additively we have an isomorphism
$H_{1}(Y_n;\bz)\simeq \bz^{n-r}\bigoplus \bz_2^{r}$ where $r$ is the
number of $n+1\leq q\leq 2n$ for which there exists a $p<q$ with
$c_{p-n,q-n}=1$. \epropo

\bcor The commutator subgroup $[\pi_1(Y_n),\pi_1(Y_n)]$ is a free
abelian group with generators $\alpha_q^2$ where $n+1\leq q\leq 2n$ is
such that there exists $p<q$ with $c_{p-n,q-n}=1$. \ecor {\bf Proof:}
In $[W,W]$, $w\cdot[s_i,s_j]\cdot w^{-1}=[s_j,s_i]$ if the reduced
word $w$ contains either $s_i$ or $s_j$ but not both,
$w\cdot[s_i,s_j]\cdot w^{-1}=[s_i,s_j]$ otherwise.But we know that $[s_i,s_j]=(s_i\cdot s_j)^2=1$ when $j\neq i+n$\\
So,
$[W,W]=\langle w\cdot[s_i,s_j]\cdot w^{-1}\ |\ 1\leq i,j\leq d,\ w\in
W\rangle=\langle (s_i\cdot s_{i+n})^2\ |\ 1\leq i\leq n\rangle.$
Observe that
$\alpha_j^2=(s_j\cdot s_1^{c_{j-n,1}}\cdots
s_n^{c_{j-n,n}})^2=(s_j\cdot s_{j-n})^2$
for $j=n+1\leq j\leq 2n$.  Thus $[W,W]$ is generated by
$\langle \alpha_j^2\ |\ n+1 \leq j\leq 2n\rangle$ as a subgroup of
$\pi_{1}(Y_n)$. Moreover, by \cite[Lemma 4.2]{u}, the torsion
elements of $W$ and hence in $[W,W]$ are only of order $2$.  Hence it
follows that $[W,W]$ is a free abelian group generated by
$\alpha_j^2, n+1 \leq j\leq 2n$.  Since $[\pi_{1}(Y_n),\pi_{1}(Y_n)]$
is a subgroup of $[W,W]$, by Proposition \ref{presfirsthomology}, the
corollary follows.  $\Box$

\bpropo\label{nil} The group $\pi_1(Y_n)$ is nilpotent if and only if
it is abelian. \epropo {\bf Proof:} Let $\bar{\alpha_q}^2$,
$n+1\leq q\leq 2n$ is such that there exists a $p<q$ such that
$c_{p-n,q-n}=1$. Then we see that
$\alpha^2_q \alpha_p\alpha^{-2}_q\alpha_p^{-1}=(s_{q}s_{q-n})^4=\alpha_q^4$,
which belongs to $[\pi_{1}(Y_n),[\pi_1(Y_n),\pi_1(Y_n)]]$. Proceeding
similarly by induction we get that
\[\alpha^{2(k-1)}_q \alpha_p\alpha^{-2(k-1)}_q\alpha_p^{-1}=\alpha_q^{2k}\]
belongs to $\pi_{1}(Y_n)^{(k)}$. Here
$\pi_1(Y_n)^{(1)}:=[\pi_1(Y_n),\pi_1(Y_n)]$ and
$\pi_{1}(Y_n)^{(k)}:=[\pi_1(Y_n), \pi_1(Y_n)^{(k-1)}]$. Since $\alpha_q$ is
of infinite order in $W$, the proposition follows.  $\Box$

\brem\label{nilpr} In particular, by Corollary \ref{ab}, $\pi_1(Y_n)$
is nilpotent if and only if $Y_n$ is a product of $n$-copies of
$\bp^1$'s.  \erem
 
\brem Here we mention that the fundamental group of a real toric
variety is not in general solvable. For example, the non-orientable
surfaces of genus $g$ are real toric varieties (see for example
\cite[Section 4.5, Remark 4.5.2]{j} and \cite[Remark 3.3]{u}), whose
fundamental groups contain free subgroups of rank $g-1$ (see for
example \cite[p. 62, Section 4]{bmms}). Thus whenever $g\geq 3$, these
groups are not solvable. Hence the property proved in Proposition
\ref{sol} is specific to the real Bott tower.  \erem

\section{Cohomology ring with $\bz_2$-coefficients}
In this section we shall describe cohomology ring $H^*(Y_n;\bz_2)$ in
terms of generators and relations by applying \cite[Theorem 5.12]{dj},
viewing $Y_n$ as a small cover.  More precisely, we state the
following theorem.

\bth\label{cohom} Let $\cR:=\bz_2[x_1,x_2\ldots, x_{2n}]$ and let $\ci$
denote the ideal in $\cR$ generated by the following set of elements
\be\label{st}\{x_ix_{n+i}~,~x_i-x_{n+i}+\sum_{j=1}^{i-1}c_{j,i}\cdot
x_{n+j} ~~\forall~1\leq i\leq n \}.\ee As a graded $\bz_2$-algebra
$H^*(Y_n;\bz_2)$ is isomorphic to $\cR/\ci$.  \eeth

{\bf Proof:} The proof immediately follows by applying 
\cite[Theorem 4.13]{dj} viewing $Y_n$ as a small cover over an $n$-cube with the
characteristic map given by means of the Bott matrix.  Thus we get the
following isomorphism of graded $\bz_2$-algebras:
\be\label{mainiso}\psi:\cR/\ci\simeq H^*(Y_n;\bz_2) \ee
where $\psi$ maps $x_i+I$ to the fundamental class of the {\it
  characteristic submanifold} $[M_{i}]\in H^{1}(Y_n;\bz_2)$ for
$1\leq i\leq 2n$. Hence the theorem.$\Box$

\bcor\label{simcohom}
Let $\wt{\cR}:=\bz_2[y_1,y_2\ldots, y_{n}]$ and let $\wt{\ci}$
denote the ideal in $\wt{\cR}$ generated by the following set of elements
\begin{equation}\{y^2_i-\sum_{j=1}^{i-1}c_{j,i}\cdot y_{i}y_{j}
  ~~~\forall~1\leq i\leq n \}.\end{equation} We have the following
isomorphism of graded $\bz_2$-algebras:
\begin{equation}\label{gaiso}\wt{\psi}:\wt{\cR}/\wt{\ci}\simeq H^*(Y_n;\bz_2)\end{equation}
where $\wt{\psi}$ maps $y_i+I$ to the fundamental class of the {\it
  characteristic submanifold} $[M_{i+n}]\in H^{1}(Y_n;\bz_2)$ for
$1\leq i\leq n$. \ecor {\bf Proof:} In (\ref{st}), by substituting the
second relation in the first, we can simplify the presentation of
$\cR/\ci$ by reducing the generators to $x_{n+1},\ldots,x_{2n}$ and
the relations to \be\label{simprels} x_{n+i}\cdot
[x_{n+i}-\sum_{j=1}^{i-1}c_{j,i}\cdot x_{n+j}] ~~\forall~1\leq i\leq
n.\ee Now, by sending $y_i$ to $x_{i+n}$ we get the required
isomorphism from $\cR/\ci$ to $\wt{\cR}/\wt{\ci}$

\brem We can also apply results in \cite[Section 4.3]{j} to get Theorem
\ref{cohom} by viewing $Y_n$ as a real toric variety.  \erem

\brem\label{lbs} Note that, in fact there exists a canonical line
bundle $[L_{i}]$ on $Y_n$ such that $[M_i]=w_1(L_i)$ for every $1\leq
i\leq 2n$, where $w_1(L_i)$ denotes the first Stiefel Whitney class of
$L_i$. Thus the cohomology ring is generated by $w_1(L_i)$ for $1\leq
i\leq 2n$ (see \cite[Section 6.1]{dj}).\erem

\section{Stiefel Whitney classes of $Y_n$}
Let $w_{k}(Y_n)$ denote the $k$th Stiefel-Whitney class of $Y_n$ for
$0\leq k\leq n$ with the understanding that $w_0(Y_n)=1$.  Then
$w(Y_n)=1+w_1(Y_n)+\cdots+w_{n}(Y_n)$ is the total Stiefel-Whitney
class of $Y_n$.

\bth\label{rf} (i) Under the isomorphism (\ref{mainiso}) of
$H^*(Y_n;\bz_2)$ with $\cR/\ci$ we have the identification
\be\label{sw} w(Y_n)=\prod_{i=1}^{2n} (1+ x_i)\ee where $x_i$ for
$1\leq i\leq 2n$ satisfy (\ref{st}).

 (ii) We further have the following recursive formula
 \be\label{rf1}w(Y_n)=w(Y_{n-1})\cdot (1+x_n)(1+x_{2n}),\ee where
  \be\label{rels} x_n\cdot x_{2n}=0,
  x_n=x_{2n}-\sum_{i=1}^{n-1} c_{i,n}x_{n+i}.\ee 
\eeth

{\bf Proof:} The proof of (i) follows readily by applying
\cite[Corollary 6.8]{dj} for $Y_n$.

Now, we shall prove (ii). 

Note that the defining Bott matrix for $Y_{n-1}$ is the
$n-1\times n-1$ submatrix of $C$ obtained by deleting the $n$th
row and the $n$th column of $C$.

Moreover, let $\pi_n:Y_n\ra Y_{n-1}$ denote the canonical projection
of the $\br\bp^1$-bundle. Then via pullback along $\pi_n^*$,
$H^*(Y_{n-1};\bz_2)$ can be identified with the subring
$\cR'/\ci'$ of $\cR/\ci$ where
$\cR'=\bz_2[x_1,x_2,\ldots,x_{n-1},x_{n+1},\ldots, x_{2n-1} ] $ and
$\ci' $ is the ideal generated by the relations \be\label{st1}
\{x_ix_{n+i}~,~ x_i-x_{n+i}+\sum_{j=1}^{i-1}c_{j,i}\cdot x_{n+j}
~~\mbox{for}~~ 1\leq i\leq n-1\}.\ee Since $Y_n$ is an
$\br\bp^1$-bundle over $Y_{n-1}$, we further have the following
presentation of $H^*(Y_n;\bz_2)$ as an algebra over
$H^*(Y_{n-1};\bz_2)$: \be\label{presind}H^*(Y_n;\bz_2)\simeq
H^*(Y_{n-1};\bz_2)[x_n,x_{2n}]/ J\ee where $J$ is the ideal generated
by the relations \be\label{relnind} x_n\cdot x_{2n},
~~x_{n}-x_{2n}+\sum_{i=1}^{n-1}c_{i,n}x_{n+i}.\ee Furthermore, via
$\pi_n^*$ we can identify $w(Y_{n-1})$ with the
expression\be\label{sw1} w(Y_{n-1})=\prod_{i=1}^{n-1} (1+
x_i)\cdot\prod_{i=n+1}^{2n-1} (1+ x_i)\ee in $\wt{R}$ where $x_i$ for
$1\leq i\leq n-1$ and $n+1\leq i\leq 2n-1$ satisfy the relations
(\ref{st1}). Now by (\ref{sw}) and (\ref{sw1}), (ii) follows.  $\Box$

\bcor\label{simrf} \begin{enumerate}\item[(i)] The following hold in the $\bz_2$-algebra
$\cR/\ci$: \be\label{rf1'}w(Y_n)=w(Y_{n-1})(1+
\sum_{i=1}^{n-1}c_{i,n}x_{n+i}),\ee \be\label{rf2}
w_{k}(Y_n)=w_{k}(Y_{n-1})+w_{k-1}(Y_{n-1})\cdot
(\sum_{i=1}^{n-1}c_{i,n}x_{n+i}) \ee for $n\geq 2$ and
$1\leq k\leq n$. \item[(ii)] For every $1\leq k\leq n$, $w_{k}(Y_n)$ is a
$\bz_2$-linear combination of square free monomials of degree $k$ in
the variables $x_{n+1},\ldots,x_{2n-1}$ modulo $\ci$. \end{enumerate} \ecor {\bf
  Proof:} The equation (\ref{sw}) reduces to (\ref{rf1'}) by applying
(\ref{rels}). Note that under the isomorphism $\psi$ of graded
algebras $H^*(Y_n;\bz_2)$ and $\cR/\ci$,
$w_k(Y_n)\in H^k(Y_n;\bz_2)$ corresponds to a polynomial of degree $k$
in $x_i ,1\leq i\leq 2n$ modulo $\ci$ for $1\leq k\leq n$.  Thus we
get (\ref{rf2}) by comparing the degree $k$-terms on either side of
(\ref{rf1'}) and (i) follows. 

Observe that by applying (\ref{st1}), in $\cR'/\ci'$ and hence
in $\cR/\ci$, we can substitute for $x_i$ in terms of
$x_{n+1},\ldots, x_{n+i}$ modulo $\ci$ using the equality
\be\label{redvar} x_{i}=x_{n+i}+\sum_{j=1}^{i-1}c_{j,i}\cdot
x_{n+j}.\ee In particular, $w_{k}(Y_{n-1})$ (resp. $w_{k-1}(Y_{n-1})$)
can be written as a polynomial of degree $k$ (resp. $k-1$) in
$x_{n+i}, ~~1\leq i\leq n-1$.  Furthermore, multiplying either side of
(\ref{redvar}) with $x_{n+i}$, along with the equality
$x_i\cdot x_{n+i}=0$ gives \be\label{sqfree}
x_{n+i}^2=\sum_{j=1}^{i-1}c_{j,i}\cdot x_{n+i}\cdot x_{n+j}\ee for
$1\leq i\leq n-1$.  It follows that $w_{k}(Y_{n-1})$
(resp. $w_{k-1}(Y_{n-1})$) can be expressed as $\bz_2$-linear
combinations of square free monomials of degree $k$ (resp. $k-1$) in
$x_{n+i}, ~~1\leq i\leq n-1$ in the algebra $\cR/\ci$.  Now,
assertion (ii) follows readily by applying (\ref{sqfree}) again in
(\ref{rf2}).$\Box$

\bcor\label{n-1th sfclass} We have the following elegant formula for
$w_{n-1}(Y_n)$ in $H^{*}(Y_n;\bz_2)$ in terms of the Bott numbers
$c_{i,j}$:
\[ w_{n-1}(Y_n)=c_{1,2}\cdot c_{2,3}\cdots c_{n-1,n} \cdot x_{n+1}\cdot
x_{n+2}\cdots x_{2n-1}.\]
\ecor {\bf Proof:} The proof follows by induction on $n$ and
(\ref{rf2}) using the fact that in $H^*(Y_n;\bz_2)$ the folowing
relations hold : \be\label{sqrelns} x_{n+1}^2=0 ; x_{n+1}\cdot
x_{n+2}^2=0; x_{n+1}\cdot x_{n+2}\cdot x_{n+3}^2=0 \cdots
;x_{n+1}\cdot x_{n+2}\cdots x^2_{2n-2}=0 .\ee $\Box$

\subsection{Orientability of the real Bott tower}
In this section we give a necessary and sufficient condition for $Y_n$
to be orientable.

\blem\label{fswc} We have the following expression for the $w_1(Y_n)$
in $\cR/\ci$: \be\label{sw2}w_1(Y_n)=\sum_{i=1}^{n-1}(\sum_{j=i+1}^n c_{i,j})\cdot
x_{n+i}.\ee
\elem
{\bf Proof:} The lemma follows by putting $k=1$ in (\ref{rf2}) followed by
induction on $n$. 
$\Box$

\bth\label{crorient1} The space $Y_n$ is orientable if and only if
\be\label{crorient}\sum_{j=i+1}^n c_{i,j}\equiv 0\,(\bmod\, \bz_2)~~ \mbox{for every} ~~1\leq
i\leq n-1,\ee where $c_{i,j}$ are the entries of the defining Bott
matrix $C$ (see \ref{BM}). \eeth {\bf Proof:} Note that via the
isomorphism (\ref{mainiso}) of graded algebras $\cR/\ci$ and
$H^*(Y_n;\bz_2)$ $\{x_{n+1},x_{n+2},\ldots,x_{2n}\}$ corresponds to a
basis over $\bz_2$ of $H^1(Y_n;\bz_2)$. Thus by (\ref{sw2}), it
follows that $w_1(Y_n)=0$ if and only if
\be\label{crorient2}\sum_{j=i+1}^n c_{i,j}\equiv 0 \,(\bmod\, \bz_2)~~ \mbox{for every} ~~1\leq
i\leq n-1.\ee Furthermore, since a necessary and sufficient condition
for a compact connected differentiable manifold $M$ to be orientable
is $w_1(M)=0$, the theorem follows. $\Box$

\bcor\label{conseq} Let $Y_n$ be an oriented real Bott tower
$Y_n$. Then $w_{n-1}(Y_n)=0$.  \ecor {\bf Proof:} By (\ref{crorient})
it follows that $c_{n-1,n}=0$ if $Y_{n}$ is orientable. Now by
(\ref{n-1th sfclass}) the corollary follows. $\Box$

\brem The assertion of Corollary \ref{conseq} is true for any even
dimensional manifold but not in general true when the dimension is odd
(see \cite[Theorem II and examples on p. 94]{mas}). Our assertion
although specific to the case of a real Bott tower, holds in all
dimensions. \erem

\brem\label{spin} In particular, a $3$-step oriented real Bott tower
$Y_3$ satisfies $w_2(Y_3)=0$, and hence admits spin structure.  This
is a special case of the well known more general result of Steenrod
that an oriented threefold is parallelizable.\erem

\brem\label{paralell} Trivially product of $n$ copies of $\bp^1$'s is
a parallelizable Bott tower for any $n$. Converse is not true as can
be seen by the $3$-step Bott tower associated with Bott numbers
$c_{1,2}=1$, $c_{1,3}=1$ and $c_{2,3}=0$, which is parallelizable but
not a product of $\bp^1$'s.  \erem

Now we give a combinatorial characterization for $Y_n$ to admit a spin
structure.  \bth\label{crspin} The orientable Bott tower $Y_n$ admits
a spin structure if and only if in addition to (\ref{crorient}) the
following identities hold for $1\leq j<k\leq n-2$: \be\label{crspin1}
‎‎\sum_{r=j+1}^ n\sum _{\substack {s=k+1\\ s\neq r}}^n c_{j,r}\cdot
c_{k,s}+c_{j,k}‎‎\sum_{\substack{r,s=k+1 \\ r< s}}^n c_{k,r}\cdot
c_{k,s}\equiv 0\,(\bmod\, \bz_2)\ee where $c_{i,j}$ are as defined in
(\ref{BM}).  \eeth {\bf Proof:} By Theorem \ref{cohom} and by equation
(\ref{sw}), $w(Y_n)$ can be identified with the class in $\cR/\ci$ of
the following term \be\label{tsw} \prod_{j=1}^n \left(1+ x_j+ x_{n+j}+
  x_j\cdot x_{n+j}\right) .\ee Further, using the relations (\ref{st})
in $\ci$ we can rewrite (\ref{tsw}) as \be\label{split} \prod_{j=2}^n
\left(1+\sum_{i=1}^{j-1} c_{i,j}\cdot x_{n+i}\right) \ee Furthermore,
since Theorem \ref{cohom} gives an isomorphism of graded
$\bz_2$-algebras, the degree $2$ term of $w(Y_n)$, namely $w_2(Y_n)$
can be identified with the degree $2$ term of expression (\ref{split})
which is the class of the following term in $ \cR/\ci$:
\be\label{fsw2}\sum_{1\leq j<k\leq n-1} \left(\sum_{r=j+1}^
  n\sum_{\substack{s= k+1 \\ s\neq r}}^ n c_{j,r} c_{k,s} ‎\right)
x_{n+j}x_{n+k}‎+\sum_{k=1}^{n-2}\left(\sum_{\substack{r,s=k+1 \\ r<
      s}}^n c_{k,r} c_{k,s}\right) x_{n+k}^2.\ee By substituting
$x_{n+k}^2=\sum_{j=1}^{k-1}c_{j,k} \cdot x_{n+j} \cdot x_{n+k}$ from
(\ref{sqfree}) and $c_{n-1,n}=0$ from (\ref{crorient}) in
(\ref{fsw2}), we get that $w_2(Y_n)$ can be identified with the class
of the following term in $\cR/\ci$ \be\label{fsw2'}\sum_{1\leq j<k\leq
  n-2} \left(\sum_{r=j+1}^ n\sum_{\substack {s= k+1 \\ s\neq r}}^ n
  c_{j,r} c_{k,s}
  ‎+c_{j,k}\cdot \sum_{\substack{r,s=k+1 \\
      r< s}}^n c_{k,r} c_{k,s}\right) x_{n+j}x_{n+k}.\ee Further, as a
graded $\bz_2$-vector space $H^2(Y_n;\bz_2)$ is isomorphic to the
subspace of $\cR/\ci$ freely generated over $\bz_2$ by the classes of
$x_{n+j}x_{n+k}$, $1\leq j<k\leq n$. Moreover, the necessary and
sufficient condition for an orientable manifold $M$ to admit a spin
structure is $w_2(M)=0$.  Thus the theorem follows from (\ref{fsw2'}).
$\Box$ \beg\label{4-step with spin} The $4$-step Bott towers admitting
spin structure are classified by the following list of associated Bott
matrices.
\begin{equation*}
  \begin{pmatrix}
   1 & 0 & 0 & 0\\
   0 & 1 & 0 & 0\\
   0 & 0 & 1 & 0\\
   0 & 0 & 0 & 1
   \end{pmatrix}
   \quad
\begin{pmatrix}
   1 & 0 & 0 & 0\\
   0 & 1 & 1 & 1\\
   0 & 0 & 1 & 0\\
   0 & 0 & 0 & 1
   \end{pmatrix}
 \quad
\begin{pmatrix}
   1 & 0 & 1 & 1\\
   0 & 1 & 0 & 0\\
   0 & 0 & 1 & 0\\
   0 & 0 & 0 & 1
   \end{pmatrix}
   \quad
   \begin{pmatrix}
   1 & 0 & 1 & 1\\
   0 & 1 & 1 & 1\\
   0 & 0 & 1 & 0\\
   0 & 0 & 0 & 1
   \end{pmatrix}
 \end{equation*}
 
 \begin{equation*}
  \begin{pmatrix}
   1 & 1 & 1 & 0\\
   0 & 1 & 0 & 0\\
   0 & 0 & 1 & 0\\
   0 & 0 & 0 & 1
   \end{pmatrix}
 \quad
\begin{pmatrix}
   1 & 1 & 1 & 0\\
   0 & 1 & 1 & 1\\
   0 & 0 & 1 & 0\\
   0 & 0 & 0 & 1
   \end{pmatrix}
   \quad
   \begin{pmatrix}
   1 & 1 & 0 & 1\\
   0 & 1 & 0 & 0\\
   0 & 0 & 1 & 0\\
   0 & 0 & 0 & 1
   \end{pmatrix}
 \quad
  \begin{pmatrix}
   1 & 1 & 0 & 1\\
   0 & 1 & 1 & 1\\
   0 & 0 & 1 & 0\\
   0 & 0 & 0 & 1
   \end{pmatrix}
 \end{equation*}
 \eeg \brem Note that the above list of Bott matrices exhausts all
 orientable $4$-step Bott towers. Thus it follows that every
 orientable $4$-step Bott tower is also spin. Moreover, it is known
 that a $4$-manifold $M$ is parallelizable if and only if it admits a
 spin structure (i.e $w_1(M)=w_2(M)=0$) and has vanishing Euler
 characteristic and signature ($\chi(M)=\sigma(M)=0$) (see
 \cite[Section 4]{hh} and \cite[p. 699]{k}). Moreover, by Hirzebruch
 signature formula, $\sigma(M)=\frac{1}{3}p_1(M)[M]$, where $p_1(M)$
 is the first Pontrjagin class and $[M]$ the fundamental class of
 $M$. Now, a real Bott tower has vanishing Euler characteristic (see
 Remark \ref{ec}) and vanishing Pontrjagin classes by (\cite[Corollary
 6.8 (i)]{dj}). Thus it follows that a $4$-step Bott tower is
 orientable if and only if it is parallelizable. Further, it
 corresponds to one of the eight Bott matrices in the above
 list. \erem
 
The following example shows that this is not the case in dimensions
$5$ and higher. Indeed there are $n$-step Bott towers which are
orientable but not spin when $n\geq 5$.

\beg Let $Y_n$ be the $n$-step Bott tower, $n\geq 5$, associated to
the Bott numbers $c_{1,2}=1$, $c_{1,n-2}=1$; $c_{n-2,n-1}=1$,
$c_{n-2,n}=1$; $c_{i,j}=0$ for $i\neq 1,n-2$ and $c_{1,j}=0$ for
$j\neq 2,n-2$.  These numbers clearly satisfy (\ref{crorient}) but not
(\ref{crspin1}).  Indeed in this case, when $j=1$ and $k=n-2$, the left
hand side of (\ref{crspin1}) is \\
$c_{1,2}\, (c_{n-2,n-1}+c_{n-2,n})+c_{1,n-2}\,
(c_{n-2,n-1}+c_{n-2,n})+c_{1,n-2}\, c_{n-2,n-1}\,
c_{n-2,n} \equiv 1\,(\bmod\, \bz_2).$
\eeg
\bdefe
We call the Bott matrix $C$ {\em spin} if and only if the associated Bott
tower $Y_n=Y(C)$ is spin.
\edefe

Let $R_i$ denote the $i$th row vector
$(0,\ldots,0,1=c_{i,i},c_{i,i+1}, c_{i,i+2},\ldots, c_{i,n})$ of
$C$. For every $1\leq j<k\leq n$, we define an $n\times n$ Bott matrix
$C_{jk}$ with $R_j$ as the $j$th row and $R_k$ as the $k$th
row. Further, for $i\neq j,k$, we let the $i$th row of $C_{jk}$ have
$1$ as the $(i,i)$th entry and all other entries as $0$.

\bcor\label{crspinsub} The Bott matrix $C$ is spin if and only if
$C_{jk}$ is spin for every $1\leq j<k\leq n-2$.  \ecor {\bf Proof:}
From Theorem \ref{crspin} a necessary and sufficient condition for $C$
to be spin is that the entries $c_{i,j}, ~i+1\leq j\leq n$ on the row
$R_i$ for every $1\leq i\leq n$ satisfy (\ref{crorient}) and further,
the entries $c_{j,r}$, $j+1\leq r\leq n$ of $R_j$ and $c_{k,s}$,
$k+1\leq s\leq n$ of $R_k$ for every $1\leq j<k\leq n$ satisfy
(\ref{crspin1}).

Again by Theorem \ref{crspin} it follows that, the necessary and
sufficient condition for the Bott matrix $C_{jk}$ to be spin is that
the entries $c_{j,r}$, $j+1\leq r\leq n$ of the $j$th row and the
entries $c_{k,s}$, $k+1\leq s\leq n$ of the $k$th row of $C_{jk}$,
satisfy (\ref{crorient}) and (\ref{crspin1}).  This can be readily
seen because any row of $C_{jk}$ other than the $j$th or $k$th row has
all entries as $0$ except the diagonal entry which is $1$. Thus the
entries above the main diagonal on the $i$th row of $C_{jk}$ where
$i\neq j,k$, trivially satisfy (\ref{crorient}). Moreover, if either
$i\neq j,k$ or $l\neq j,k$ and $1\leq i<l\leq n$, the entries of
$C_{jk}$ above the main diagonal, on the $i$th and the $l$th row
trivially satisfy (\ref{crspin1}). Hence the corollary. $\Box$

\brem\label{compwithg} See \cite{g} for results on criterion for $Y_n$ to admit spin
structure using different methods. We however note here that the main
result \cite[Theorem 1.2]{g} follows immediately from Corollary
\ref{crspinsub} above. This can be seen because the conditions
(\ref{crorient}) and (\ref{crspin1}) need to be checked only for the
entries of pairs of nonzero rows of the matrix $A$ where $A:=C-I$. In
\cite{g}, $A$ is called the Bott matrix. Thus we need to check the
spin condition only on those $A_{jk}:=C_{jk}-I$ with $j$th and $k$th
nonzero rows.

Moreover, our result is more general as we do not require the
assumption in \cite[Theorem 1.2]{g} that the number of nonzero rows of
the matrix $A$ is even (see \cite[Remark 2.1]{g}).\erem We illustrate
Corollary \ref{crspinsub} by the following examples. 
\beg \begin{enumerate}\item\begin{equation*} C= \begin{pmatrix}
      1 & 0 & 0 & 0 & 0 & 0 \\
      0 & 1 & 0 & 0 & 1 & 1 \\
      0 & 0 & 1 & 0 & 1 & 1 \\
      0 & 0 & 0 & 1 & 1 & 1 \\
      0 & 0 & 0 & 0 & 1 & 0 \\
      0 & 0 & 0 & 0 & 0 & 1
           \end{pmatrix} 
\end{equation*}
\begin{equation*}       
\underbrace{\begin{pmatrix}
            1 & 0 & 0 & 0 & 0 & 0 \\
            0 & 1 & 0 & 0 & 1 & 1 \\
            0 & 0 & 1 & 0 & 1 & 1 \\
            0 & 0 & 0 & 1 & 0 & 0 \\
            0 & 0 & 0 & 0 & 1 & 0 \\
            0 & 0 & 0 & 0 & 0 & 1
           \end{pmatrix}}_{C_{23}} 
\quad    
\underbrace{\begin{pmatrix}
            1 & 0 & 0 & 0 & 0 & 0 \\
            0 & 1 & 0 & 0 & 1 & 1 \\
            0 & 0 & 1 & 0 & 0 & 0 \\
            0 & 0 & 0 & 1 & 1 & 1 \\
            0 & 0 & 0 & 0 & 1 & 0 \\
            0 & 0 & 0 & 0 & 0 & 1
           \end{pmatrix}}_{C_{24}}         
\quad 
\underbrace{\begin{pmatrix}
             1 & 0 & 0 & 0 & 0 & 0 \\
             0 & 1 & 0 & 0 & 0 & 0 \\
             0 & 0 & 1 & 0 & 1 & 1 \\
             0 & 0 & 0 & 1 & 1 & 1 \\
             0 & 0 & 0 & 0 & 1 & 0 \\
             0 & 0 & 0 & 0 & 0 & 1
           \end{pmatrix}}_{C_{34}}
  \end{equation*}
  Clearly $\sum_{j=i+1}^n c_{i,j}\equiv 0\pmod {\bz_2}$ for all $1\leq i\leq 5$. \\
  When $j=2$ and $k=3$ the left hand side of (\ref{crspin}) is\\ $c_{2,3}(c_{3,4}+c_{3,5}+c_{3,6}) +
  c_{2,4}(c_{3,5}+c_{3,6})+ c_{2,5}(c_{3,4}+c_{3,6}) + c_{2,6}(c_{3,4}+c_{3,5})
  +\\ c_{2,3}(c_{3,4}\,c_{3,5}+c_{3,4}\,c_{3,6}+c_{3,5}\,c_{3,6}) \equiv
  0\,(\bmod\, \bz_2)$

  When $j=3$ and $k=4$ the left hand side of (\ref{crspin}) is \\
  $c_{3,4}(c_{4,5}+c_{4,6})+ c_{3,5}\,c_{4,6}+c_{3,6}\,c_{4,5} +c_{3,4}\cdot
  c_{4,5}\cdot c_{4,6}\equiv 0\,(\bmod\, \bz_2)$
 
When $j=2$ and $k=4$ the left hand side of (\ref{crspin}) is \\
$c_{2,3}(c_{4,5}+c_{4,6})+c_{2,4}(c_{4,5}+c_{4,6})+c_{2,5}\,c_{4,6}+c_{2,6}\,c_{4,5}
+c_{2,4}\,c_{4,5}\,c_{4,6} \equiv 0\\\,(\bmod\, \bz_2)$.

Since $C_{23}, C_{24},C_{34}$ are all spin by Corollary
\ref{crspinsub}, $C$ is spin. We do not consider the matrices $C_{1l}$
for $2\leq l\leq 4$, since the first row of $C$ has nonzero entry only
on the main diagonal.
\item \begin{equation*} C= \begin{pmatrix}
      1 & 0 & 1 & 1 & 0 \\
      0 & 1 & 1 & 1 & 0 \\
      0 & 0 & 1 & 1 & 1 \\
      0 & 0 & 0 & 1 & 0 \\
      0 & 0 & 0 & 0 & 1
           \end{pmatrix} 
           \end{equation*}
        \begin{equation*}
          \underbrace{\begin{pmatrix}
            1 & 0 & 1 & 1 & 0 \\
            0 & 1 & 1 & 1 & 0 \\
            0 & 0 & 1 & 0 & 0 \\
            0 & 0 & 0 & 1 & 0 \\
            0 & 0 & 0 & 0 & 1
                \end{pmatrix}}_{C_{12}}
\quad        
   \underbrace{\begin{pmatrix}
                1 & 0 & 1 & 1 & 0 \\
                0 & 1 & 0 & 0 & 0 \\
                0 & 0 & 1 & 1 & 1 \\
                0 & 0 & 0 & 1 & 0 \\
                0 & 0 & 0 & 0 & 1
                \end{pmatrix}}_{C_{13}}             
\quad
   \underbrace{\begin{pmatrix}
                1 & 0 & 0 & 0 & 0 \\
                0 & 1 & 1 & 1 & 0 \\
                0 & 0 & 1 & 1 & 1 \\
                0 & 0 & 0 & 1 & 0 \\
                0 & 0 & 0 & 0 & 1
                \end{pmatrix}}_{C_{23}}     
       \end{equation*}
       Clearly $\sum_{j=i+1}^n c_{i,j}=0$ for all $1\leq i\leq 4$.\\
       When $j=1$ and $k=2$ the left hand side of (\ref{crspin}) is\\
       $c_{1,2}(c_{2,3}+c_{2,4}+c_{2,5})+c_{1,3}(c_{2,4}+c_{2,5})+c_{1,4}(c_{2,3}+c_{2,5})+c_{1,5}(c_{2,3}+c_{2,4})
       +\\ c_{1,2}(c_{2,3}\,c_{2,4}+c_{2,3}\,c_{2,5}+c_{2,4}\,c_{2,5})
       \equiv 0 \,(\bmod\, \bz_2)$.\\
       When $j=1$ and $k=3$ the left hand side of (\ref{crspin}) is\\
       $c_{1,2}(c_{3,4}+c_{3,5})+c_{1,3}(c_{3,4}+c_{2,5})+c_{1,4}\,c_{3,5}+
       c_{1,5}\,c_{3,4}+
       c_{1,3}\,c_{3,4}\,c_{3,5} \equiv 1\\\,(\bmod\, \bz_2)$.
       Thus $C$ is not spin since $C_{13}$ is not spin.
\item \begin{equation*}
        C= \begin{pmatrix}
            1 & 0 & 0 & 0 & 0 & 0 & 0 \\
            0 & 1 & 1 & 1 & 1 & 1 & 0\\
            0 & 0 & 1 & 1 & 1 & 1 & 1\\
            0 & 0 & 0 & 1 & 0 & 1 & 1\\
            0 & 0 & 0 & 0 & 1 & 0 & 0\\
            0 & 0 & 0 & 0 & 0 & 1 & 0\\
            0 & 0 & 0 & 0 & 0 & 0 & 1
           \end{pmatrix} 
    \end{equation*}
      \begin{equation*}   
     \underbrace{\begin{pmatrix}
            1 & 0 & 0 & 0 & 0 & 0 & 0 \\
            0 & 1 & 1 & 1 & 1 & 1 & 0\\
            0 & 0 & 1 & 1 & 1 & 1 & 1\\
            0 & 0 & 0 & 1 & 0 & 0 & 0\\
            0 & 0 & 0 & 0 & 1 & 0 & 0\\
            0 & 0 & 0 & 0 & 0 & 1 & 0\\
            0 & 0 & 0 & 0 & 0 & 0 & 1
           \end{pmatrix} }_{C_{23}}
\quad
\underbrace{\begin{pmatrix}
            1 & 0 & 0 & 0 & 0 & 0 & 0 \\
            0 & 1 & 1 & 1 & 1 & 1 & 0\\
            0 & 0 & 1 & 0 & 0 & 0 & 0\\
            0 & 0 & 0 & 1 & 0 & 1 & 1\\
            0 & 0 & 0 & 0 & 1 & 0 & 0\\
            0 & 0 & 0 & 0 & 0 & 1 & 0\\
            0 & 0 & 0 & 0 & 0 & 0 & 1
           \end{pmatrix}}_{C_{24}} 
\quad   
\underbrace{\begin{pmatrix}
            1 & 0 & 0 & 0 & 0 & 0 & 0 \\
            0 & 1 & 0 & 0 & 0 & 0 & 0\\
            0 & 0 & 1 & 1 & 1 & 1 & 1\\
            0 & 0 & 0 & 1 & 0 & 1 & 1\\
            0 & 0 & 0 & 0 & 1 & 0 & 0\\
            0 & 0 & 0 & 0 & 0 & 1 & 0\\
            0 & 0 & 0 & 0 & 0 & 0 & 1
           \end{pmatrix}}_{C_{34}} 
       \end{equation*}\\
       Clearly $\sum_{j=i+1}^n c_{i,j}=0$ for all $1\leq i\leq 6$. \\
       When $j=2$ and $k=3$ the left hand side of (\ref{crspin}) is\\
       $c_{2,3}(c_{3,4}+c_{3,5}+c_{3,6}+c_{3,7})+c_{2,4}(c_{3,5}+c_{3,6}+c_{3,7})+
       c_{2,5}(c_{3,4}+c_{3,6}+c_{3,7})+c_{2,6}(c_{3,4}+c_{3,5}+c_{3,7} )
       +\\
       c_{2,3}(c_{3,4}\,c_{3,5}+c_{3,4}\, c_{3,6}+c_{3,4}\,
       c_{3,7}+c_{3,5}\, c_{3,6}+c_{3,5}\,c_{3,7}+c_{3,6}\, c_{3,7})
       \equiv 1\,(\bmod\, \bz_2)$. 
       Thus $C$ is not spin since $C_{23}$ is not spin.
\end{enumerate}
\eeg

\subsection{Existence of a nowhere vanishing tangent vector field}

\bpropo The $n$-step Bott tower is the total space of a fibre bundle
over $\bs^1$ with fibre an $n-1$-step Bott tower corresponding to
the Bott matrix $C^{1}$ of size $n-1\times n-1$, defined by
deleting the first row and first column of $C$. \epropo

{\bf Proof:} This can be seen as follows: Let $N'$ denote the lattice
$\bigoplus_{i=2}^{n} \bz e_i$. Let $v'_1=e_2$, $v'_2=e_3$, $\ldots$,
$v'_{n-1}=e_n$. Also let $v'_{n+1}=-e_2+c_{2,3}e_3+\cdots c_{2,n}e_n$,
$\ldots$, $v'_{2n-2}=-e_{n-1}+c_{n-1,n}e_n$ and $v'_{2n-1}=-e_n$.

We define the fan $\Delta'$ in $N'$ consisting of cones generated by
the set of vectors in any subcollection of
$\{v'_1,v'_2,\ldots,v'_{n-1},v'_{n+1},\ldots,v'_{2n-1}\}$ which does
not contain both $v'_i$ and $v'_{n+i}$ for $1\leq i\leq n-1$.

Also let $N''$ denote the lattice $\bz e_1$ and $\Delta''$ denote the
fan consisting of the one dimensional cones generated by the vectors
$e_1$ and $-e_1$ and the cone $\{0\}$ which corresponds to the real
toric variety $\bp^1_{\br}$.

By projecting to $\bz e_1$ we get an exact sequence of fans
\be\label{es} 0\ra (\Delta',N')\ra (\Delta,N)\ra (\Delta'',N'')\ra 0
.\ee Moreover, the fan $\wt{\Delta''}$ in $N$ consisting of the one
dimensional cones generated by the vectors $e_1$ and
$-e_1+c_{1,2}e_2+\cdots+c_{1,n}e_n$ and the cone $\{0\}$ is a lift of
$\Delta''$. Moreover, it can be seen that every cone $\sigma$ of
$(N,\Delta)$ is a sum $\sigma'+\sigma''$ of a cone in $(N',\Delta')$
and a cone in $(N,\wt{\Delta''})$.

Thus we see that $Y_n=X(\Delta,N)$ is a toric fibre bundle over
$\bs^1$ with fibre an $n-1$-step real Bott tower corresponding to the
fan $(N',\Delta')$ (see \cite[p. 41]{f}). $\Box$

In this convention, we shall denote the $n$-step Bott tower by $Z_n$
and the fibre, which is the $n-1$-step real Bott tower associated to
the matrix $C^{1}$ and hence the fan $(N',\Delta')$, by $Z_{n-1}$. We
can iterate this process and view $Z_{n-1}$ again as a fibre bundle
over $\bs^1$ with fibre $Z_{n-2}$ which is the $n-2$-step real Bott
tower associated to the matrix $C^{2}$ of size $n-2\times n-2$,
obtained by deleting the first and the second rows and columns of $C$.
Continuing this process $n-2$ times we finally get that $Z_2$ is a two
step Bott tower associated to the Bott matrix $C^{n-2}$, obtained by
deleting the first $n-2$ rows and columns of $C$. Then $Z_2$ is a fibre
bundle over $\bs^1$ with fibre $Z_1\simeq \bs^1$.

Following the above convention, in this section we shall denote the
$n$-step real Bott tower by $Z_n$. Let $p_n:Z_n\ra \bs^1$ denote the projection of this fibre bundle.

\bth\label{tvf}The $n$-step real Bott tower $Z_n$ has a nowhere
vanishing continuous tangent vector field. Moreover, the tangent
bundle splits into a direct sum of line bundles.\eeth {\bf Proof:}
Note that the tangent bundle of $Z_n$ is the direct sum of the pull
back of the tangent bundle of $\bs^1$ and the relative tangent bundle
of $Z_{n-1}$.

Since the tangent bundle of $\bs^1$ is trivial its pull back to $Z_n$
is a trivial line bundle. Thus a nowhere vanishing section of the
trivial line bundle gives a nowhere vanishing section for the tangent
bundle of $Z_n$.

Furthermore, if we assume by induction that the tangent bundle of an
$n-1$-step Bott tower is a direct sum of line bundles, then the
relative tangent bundle of $Z_{n-1}$ in the above fibration is also a
direct sum of the associated line bundles. It follows that the tangent
bundle of the $n$-step Bott tower is a direct sum of $n$-line bundles.
Hence the proof.$\Box$

\bcor\label{eulerclass} Let $Z_n$ be an orientable real Bott
tower. Then the Euler class of $Z_n$ vanishes.  \ecor {\bf Proof:} The
proof is an immediate consequence of the \cite[Theorem \ref{tvf} and
p. 101]{ms}. $\Box$

\bcor The $n$-step real Bott tower $Z_n$ is orientable (respectively
spin) implies that the successive fibres
$Z_{n-1}, Z_{n-2},\cdots, Z_{2}$ in the above iterated construction
are all orientable (respectively spin).  \ecor {\bf Proof:} This
follows from (\ref{crorient}) and (\ref{crspin1}) since the Bott
matrix corresponding to $Z_{k}$ is $C^{n-k}$ which is the matrix
obtained from $C$ by deleting the first $k$ rows and $k$
columns. $\Box$

\brem\label{ec} Alternately the assertion of Corollary \ref{eulerclass} follows
more directly from the fact that the Euler characteristic is
multiplicative for the total space of fibre bundles. Indeed
$\chi(Z_n)=0$ since $\chi(\bs^1)=0$. In fact this also implies the
statement on existence of nowhere vanishing vector field in Theorem
\ref{tvf} by Hopf's theorem.  \erem

\subsection{Real Bott-towers bound}
  
\bdefe\label{nc} A smooth compact $n$-dimensional manifold $M$ without
boundary is {\em null cobordant} if it is diffeomorphic to the
boundary of some compact smooth $n
+ 1$-dimensional manifold $\mathcal{W}$ with boundary.  \edefe

Let $w_{k}:=w_{k}(Y_n)$ for $1\leq k\leq n$. Also let $\mu_{Y_n}$
denote the fundamental class of $Y_n$ in $H_n(Y_n;\bz_2)$. Then
\be\label{swnos} \langle w_1^{r_1} \cdots w_n^{r_n},\mu_{Y_n}\rangle
\in \bz_2 \ee such that $\sum_{i=1}^n i\cdot r_i=n$ are the
Steifel-Whitney numbers of $Y_n$.

\bth\label{vanishswnos} The $n$-step Bott tower $Y_n$ is
null-cobordant.  \eeth {\bf Proof:} From Corollary \ref{simrf} (ii), it
follows that, any monomial $w_1^{r_1}\cdots w_n^{r_n}$ of total
dimension $n$, under the isomorphism $\psi$ corresponds in $\cR/\ci$,
to a $\bz_2$-linear combination of square free monomials of degree $n$
in $x_{n+1},x_{n+2},\ldots, x_{2n-1}$. But there are no square free
monomials of degree $n$ in $x_{n+j}, ~1\leq j\leq n-1$.  Thus the
monomial $w_1^{r_1}\cdots w_n^{r_n}=0$ in $H^n(Y_n;\bz_2)$ so that the
associated Stiefel-Whitney number is zero. Therefore by Thom's theorem
it follows that $Y_n$ is null-cobordant. $\Box$

\bdefe\label{onc} A smooth compact $n$-dimensional manifold $M'$
without boundary is {\em orientedly null cobordant} if it is
diffeomorphic to the boundary of some compact smooth
$n + 1$-dimensional oriented manifold $\mathcal{W'}$ with boundary.\edefe

Let $Y_{n}$ denote an oriented $n$-step real Bott tower.  Let
$p_{i}:=p_i(Y_{n})$ denote the $i$th Pontrjagin class of $Y_{n}$ in
$H^{4i}(Y_{n},\bz)$ and $\mu_{Y_{n}}$ denote the fundamental homology
class in $H_{n}(Y_{n},\bz)$. Then for each $I=i_1,\ldots,i_r$ a
partition of $k$, the $I$th Pontrjagin number of $Y_{n}$ is given by
\be\label{pnos} \langle p_{i_1}\cdots p_{i_r},\mu_{Y_{n}}\rangle\in
\bz\ee when $n=4k$. It is zero when $n$ is not divisible by $4$ (see
\cite[p. 185]{ms}).  \bcor\label{orcobor} Let $Y_{n}$ be an oriented
$n$-step real Bott tower, then it is orientedly null-cobordant.  \ecor
{\bf Proof:} Note that \cite[Corollary 6.8 (i)]{dj}, implies that all
the Pontrjagin numbers of $Y_{n}$ vanish. Moreover, we have shown
above in the proof of Theorem \ref{vanishswnos} that all the
Stiefel-Whitney numbers of $Y_{n}$ vanish. Thus the corollary follows
by Wall's theorem (\cite[Section 8, Corollary 1]{wall}). $\Box$

\brem There are examples of real toric varieties whose top
Stiefel-Whitney class does not vanish. For example the non-orientable
surfaces of odd genus are real toric varieties (see \cite[Section 4.5,
Remark 4.5.2]{j}
\cite[Remark 3.3]{u}) having non-vanishing second Stiefel-Whitney class.  Thus the
properties we have proved in this and the preceding section are all
specific to real Bott towers.\erem

\brem\label{rbsfm} Another motivation for this work was to relate the
topology of the real Bott tower with the topology of the real Bott
Samelson manifolds using the degeneration results of Grossberg and
Karshon (see \cite{gk}). This could further be used to obtain results
on the topology of real flag manifolds. In particular, to study the
fundamental group and cohomology ring of the real Bott Samelson
manifolds and real flag variety. This shall again be taken up in
future work.  \erem

\brem During the process of this work the authors also came across the
work of Kamishima and Masuda \cite{km} where among other results they
also compute the fundamental group and cohomology ring of real Bott
towers but using different methods. \erem

\section*{Appendix}

Proof that $x_{p,q}^\epsilon$ is a conjugate of either $x_{p,q}$ or $x_{p,q}^{-1}$ by an element of the free group generated by $S$ :\\
\begin{enumerate}
 \item {$c_{p-n,q-n}=0$, $\epsilon_{p-n}=0$ and $\epsilon_{q-n}=0$} :
 \begin{equation*}
  \begin{split}
   x_{p,q}^\epsilon & =\alpha_p\,\alpha_q\,\alpha_p^{-1}\,\alpha_q^{-1} \\
   ~ & =x_{p,q}
  \end{split}
 \end{equation*}
 \item {$c_{p-n,q-n}=0$, $\epsilon_{p-n}=0$ and $\epsilon_{q-n}=1$} :
 \begin{equation*}
  \begin{split}
   x_{p,q}^\epsilon & = \alpha_p\,\alpha_q^{-1}\,\alpha_p^{-1}\,\alpha_q\\
   ~ & = \alpha_q^{-1}\cdot(\alpha_q\,\alpha_p\,\alpha_q^{-1}\,\alpha_p^{-1})\cdot\alpha_q\\
   ~ & = \alpha_q^{-1}\cdot(\alpha_p\,\alpha_q\,\alpha_p^{-1}\,\alpha_q^{-1})^{-1}\cdot\alpha_q\\
   ~ & = \alpha_q^{-1}\cdot x_{p,q}^{-1}\cdot\alpha_q
  \end{split}
 \end{equation*}
\item {$c_{p-n,q-n}=0$, $\epsilon_{p-n}=1$ and $\epsilon_{q-n}=0$} :
 \begin{equation*}
  \begin{split}
   x_{p,q}^\epsilon & = \alpha_p^{-1}\,\alpha_q\,\alpha_p\,\alpha_q^{-1}\\
   ~ & = \alpha_p^{-1}\cdot(\alpha_q\,\alpha_p\,\alpha_q^{-1}\,\alpha_p^{-1})\cdot\alpha_p\\
   ~ & = \alpha_p^{-1}\cdot(\alpha_p\,\alpha_q\,\alpha_p^{-1}\,\alpha_q^{-1})^{-1}\cdot\alpha_p\\
   ~ & = \alpha_p^{-1}\cdot x_{p,q}^{-1}\cdot\alpha_p
  \end{split}
 \end{equation*}
\item {$c_{p-n,q-n}=0$, $\epsilon_{p-n}=1$ and $\epsilon_{q-n}=1$} :
 \begin{equation*}
  \begin{split}
   x_{p,q}^\epsilon & = \alpha_p^{-1}\,\alpha_q^{-1}\,\alpha_p\,\alpha_q\\
   ~ & = (\alpha_p\,\alpha_q)^{-1}\cdot(\alpha_p\,\alpha_q\,\alpha_p^{-1}\,\alpha_q^{-1})\cdot(\alpha_p\,\alpha_q)\\
   ~ & = (\alpha_p\,\alpha_q)^{-1}\cdot x_{p,q}\cdot(\alpha_p\,\alpha_q)
  \end{split}
 \end{equation*}
\item {$c_{p-n,q-n}=1$, $\epsilon_{p-n}=0$ and $\epsilon_{q-n}=0$} :
 \begin{equation*}
  \begin{split}
   x_{p,q}^\epsilon & =\alpha_p\,\alpha_q^{-1}\,\alpha_p^{-1}\,\alpha_q^{-1} \\
   ~ & =x_{p,q}
  \end{split}
 \end{equation*}
 
\item {$c_{p-n,q-n}=1$, $\epsilon_{p-n}=1$ and $\epsilon_{q-n}=0$} :
 \begin{equation*}
  \begin{split}
   x_{p,q}^\epsilon & = \alpha_p^{-1}\,\alpha_q^{-1}\,\alpha_p\,\alpha_q^{-1}\\
   ~ & = (\alpha_p\,\alpha_q^{-1})^{-1}\cdot(\alpha_p\,\alpha_q^{-1}\,\alpha_p^{-1}\,\alpha_q^{-1})\cdot(\alpha_p\,\alpha_q^{-1})\\
   ~ & = (\alpha_p\,\alpha_q^{-1})^{-1}\cdot x_{p,q}\cdot(\alpha_p\,\alpha_q^{-1})
  \end{split}
 \end{equation*}
\item {$c_{p-n,q-n}=1$, $\epsilon_{p-n}=1$ and $\epsilon_{q-n}=1$} :
 \begin{equation*}
  \begin{split}
    x_{p,q}^\epsilon & = \alpha_p^{-1}\,\alpha_q\,\alpha_p\,\alpha_q\\
   ~ & = \alpha_p^{-1}\cdot(\alpha_q\,\alpha_p\,\alpha_q\,\alpha_p^{-1})\cdot\alpha_p\\
   ~ & = \alpha_p^{-1}\cdot(\alpha_p\,\alpha_q^{-1}\,\alpha_p^{-1}\,\alpha_q^{-1})^{-1}\cdot\alpha_p\\
   ~ & = \alpha_p^{-1}\cdot x_{p,q}^{-1}\cdot\alpha_p
   \end{split}
 \end{equation*}
\end{enumerate}


\begin{thebibliography}{9}
\bibitem {bb} A. Bahri and M. Bendersky, {\it The $KO$-theory of toric 
manifolds,} Trans. Amer. Math. Soc. {\bf 352}, 1191--1202.  

\bibitem{bmms} K. Bencsath et al. {\it Aspects of the theory of free
    groups}, Algorithmic Problems in Groups and Semigroups edited by
  Jean-Camille Birget, Stuart Margolis, John Meakin, Mark V. Sapir, 
  Trends in Mathematics, Springer, 51-90.
    
\bibitem{bs} R. Bott and H. Samelson, {\it Applications of the theory
    of Morse to symmetric spaces}, Amer. J. Math. {\bf 80}, 1958,
  964--1029.

\bibitem{br} V. Buchstaber and N. Ray, {\it Tangential structures on
    toric manifolds and connected sums of polytopes,}
  Internat. Math. Res. Notices, {\bf 4}, (2001), 193-219.

\bibitem{cr} Y. Civan and N. Ray, {\it Homotopy decompositions and
    $K$-theory of Bott towers,} $K$-Theory, {\bf 34}, (2005), no. 1,
  1--33.

\bibitem{dj} M. W. Davis and T. Januszkiewicz, Convex polytopes,
  Coxeter orbifolds and torus actions, {\it Duke Math. J.} {\bf 62}
  (1991), 417-451.

\bibitem{du} Haibao Duan, Multiplicative rule of Schubert class, {\it
    Invent. Math.}, {\bf 159}, no.2, (2005) pp.407-436.

\bibitem{f} W. Fulton, Introduction to toric varieties, Ann Math Studies {\bf 
131},(1993), Princeton Univ. Press, Princeton, NJ. 

\bibitem{cmo} Suyoung Choi, Mikiya Masuda, Sang-il Oum, Classification
  of real Bott manifolds and acyclic digraphs math.arXiv:1006.4658.

\bibitem{cms} Suyoung Choi, Mikiya Masuda and Dong Youp Suh,
  Topological classification of generalized Bott towers {\it
    Trans. Amer. Math. Soc.}{\bf 362} (2010), 1097-1112.

\bibitem{hh} Hirzebruch, F und Hopf, H, Felder von Fl$\mbox{\"{a}}$chenelementen in
  $4$-dimensionalen Mannigfaltigkeiten, {\it Math. Annalen}, {\bf 136},
  (1958), 156-172.

\bibitem{g} A. G$\mbox{\c{a}}$sior, Spin-structures on real Bott
  manifolds, math.arXiv:1506.06884v5

\bibitem{gk} M. Grossberg and Y. Karshon, Bott towers complete
integrability and the extended character of representations, {\it Duke
Math. J.} {\bf 76} (1994) 23-58.

\bibitem{mas} W. S. Massey, On the Steifel-Whitney classes of a
  manifold, {\it Amer. J. Math.}, {\bf 82}, No 1, (1960), 92-102.

\bibitem{m1} M. Masuda, Classification of real Bott manifolds,
  math. arXiv: 0809:2178.


\bibitem{mp} M.  Masuda, T.  E.  Panov, Semifree circle actions, Bott
  towers and quasitoric manifolds, {\it Mat.  Sb.}, {\bf 199}, Number 8
  (2008), 95–122.

\bibitem{ms} J. W. Milnor and J.D. Stasheff, Characteristic Classes,
  Hindustan Book Agency, TRIM Series {\bf 32} (2005).

 \bibitem{j} J. Jurkiewicz, Torus embeddings, polyhedra, $k^*$-actions
and homology,\\ {\it Dissertationes Math.} {\bf 236} (1985) 1-69.

\bibitem{k} J. Korba$\mbox{\v{s}}$, Distributions, vector distributions, immersions
  of manifolds in Euclidean spaces, {\it Handbook of Global Analysis}
  edited by Demeter Krupka, David Saunders (2008), 665-724.
	
\bibitem{km} Y. Kamishima and M. Masuda, Cohomological rigidity of real
Bott manifolds, {\it Algebr. Geom. Topol.} {\bf 9} (2009) 2479-2502.

\bibitem{su} P. Sankaran and V. Uma, K-theory of quasi-toric manifolds
  {\it Osaka J. Math.} {\bf 44}, no. 1 (2007), 71-89.


\bibitem{u} V.Uma, On the fundamental group of real toric varieties,
{\it Proc. Ind. Acad. Sci.} (Math. Sci){\bf 114} (2003) 15-31.

\bibitem{wall} C.T.C Wall, Determination of the cobordism ring, {\it
    The Annals of Mathematics}, Second Series {\bf 72} Issue 2 (1960)
  292-311.

\bibitem{w} Matthieu Willems, $K$-th$\mbox{\'{e}}$orie
  $\mbox{\'{e}}$quivariante des tours de Bott. Application
  $\mbox{\`{a}}$ la structure multiplicative de la
  𝐾-th$\mbox{\'{e}}$orie $\mbox{\'{e}}$quivariante des
  vari$\mbox{\'{e}}$t$\mbox{\'{e}}$s de drapeaux, {\it Duke
    Math. J.} {\bf 132} no. 2, (2006), 271–309.

\end{thebibliography}
\end{document}